\documentclass[11pt,openright]{memoir}
\usepackage{setspace,framed,pifont,dingbat} 
\usepackage{soul,color,amssymb,graphicx,colortbl,xcolor,amsmath} 
\usepackage{eurosym, textcomp}
\usepackage{rotating,float}
\usepackage[framed]{ntheorem}
\DisemulatePackage{setspace}
\usepackage{setspace}
\usepackage{titletoc,titlesec} 
\usepackage[T1]{fontenc} 
\usepackage[latin1]{inputenc} 
\usepackage{microtype}
\usepackage{lmodern}
\usepackage{hyperref} 
\counterwithout{section}{chapter}

\usepackage{fourier}

\newcommand{\ux}{\underline{x}}

\newcommand{\uy}{\underline{y}}

\newcommand{\uz}{\underline{z}}

\newcommand{\ut}{\underline{t}}

\newcommand{\uom}{\underline{\omega}}

\newcommand{\vphi}{\varphi}

\newcommand{\C}{\mathcal{C}\hspace*{-0.2em}\ell}

\newcommand{\RR}{\mathbb{R}}

\newcommand{\CC}{\mathbb{C}}

\numberwithin{equation}{section}

\newtheorem{thm}{Theorem}[section]
 \newtheorem{cor}[thm]{Corollary}
 \newtheorem{lem}[thm]{Lemma}
 \newtheorem{prop}[thm]{Proposition}
 \newtheorem{defn}[thm]{Definition}
 \newtheorem{rem}[thm]{Remark}
 
 \numberwithin{equation}{section}

\begin{document}

\title{The Segal-Bargmann Transform in Clifford Analysis}






\author{Swanhild Bernstein (swanhild.bernstein@math.tu-freiberg.de) and Sandra Schufmann}


%
%
\date{TU Bergakademie Freiberg, Institute of Applied Analysis}


\maketitle

\begin{abstract}

The Segal-Bargmann transform plays an essential role in signal processing, quantum physics, infinite-dimensional analysis, function theory and further topics. The connection to signal processing is the short-time Fourier transform, which can be used to describe the Segal-Bargmann transform. The classical Segal-Bargmann transform $\mathcal{B}$ maps a square integrable function to a holomorphic function square-integrable with respect to a Gaussian identity. In signal processing terms, a signal from the position space $L_2(\RR^m,\RR)$ is mapped to the phase space of wave functions, or Fock space, $\mathcal{F}^2(\CC^m,\CC)$. We extend the classical Segal-Bargmann transform to a space of Clifford algebra-valued functions. We show how the Segal-Bargmann transform is related to the short-time Fourier transform and use this connection to demonstrate that $\mathcal{B}$ is unitary up to a constant and maps Sommen's orthonormal Clifford Hermite functions $\left\{\phi_{l,k,j}\right\}$ to an orthonormal basis of the Segal-Bargmann module $\mathcal{F}^2(\CC^m,\C_m^{\CC})$. We also lay out that the Segal-Bargmann transform can be expanded to a convergent series with a dictionary of $\mathcal{F}^2(\CC^m,\C_m^{\CC})$. In other words, we analyse the signal $f$ in one basis and reconstruct it in a basis of the Segal-Bargmann module. 







\end{abstract}

\section{Introduction}

\label{intro}

Due to the importance of the Segal-Bargmann transform, there are various generalizations into quaternion and Clifford analysis. In particular, the Bargmann-Segal transformation has been studied in the theory of slice monogenic functions \cite{ACSS}, \cite{Diki},\cite{DG}, \cite{MNQ}. Our interest doesn't lie in these theories. We are interested in the importance of the Segal-Bargmann transform in its connection to the windowed Fourier transform and time-frequency analysis.

Time-frequency analysis is an important method in signal processing, because it allows to analyse a given signal simultaneously in the time and frequency domains. A well-known tool is the short-time Fourier transform. Another closely related tool is the Segal-Bargmann transform, which is our main focus in this paper.

The classical Segal-Bargmann transform maps a square integrable function to a holomorphic function square-integrable with respect to a Gaussian identity. In signal processing terms, a signal from the position space $L_2(\RR^m,\RR)$ is mapped to the phase space of wave functions $\mathcal{F}^2(\CC^m,\CC)$. In the early 1960s, V. Bargmann and I. Segal independently investigated this space \cite{Bargmann,Segal}. While Bargmann developed a theory about the space and the corresponding transform in the finite-dimensional case, Segal focused primarily on the infinite-dimensional version of the now-called Segal-Bargmann space(s) \cite{Hall}.

The space $\mathcal{F}^2(\CC^m,\CC)$ has a wide number of applications such as in infinite-dimensional analysis and stochastic distribution theory. As early as 1932, V.~Fock introduced a more general, infinite-dimensional version of this space as a quantum states space for an unknown number of particles \cite{Fock}, which is now called \textit{Fock} space. In quantum mechanics, the reproducing kernels of the Fock spaces are the so-called coherent states. Segal and Bargmann showed that an infinite union of the spaces $\mathcal{F}^2(\CC^m,\CC)$, $m\in\mathbb{N}$, is isomorphic to a certain case of the Fock space, which is why the Segal-Bargmann spaces are sometimes also called Segal-Bargmann-Fock space(s) or only Fock space. For this work, we will stick to the notion of Segal-Bargmann space.


In signal and image processing not only scalar-valued but also quaternion- and Clifford-valued signals are of interest. A monogenic signal \cite{FS,BBRH,BHRHSS}, for example, consists of a scalar-valued signal and vector components, which are the Riesz transformations of the scalar-valued signal. Other applications deal with colour images of which the colours are separated and considered as components of a Clifford-valued signal, see for example    \cite{BS,DCF,ES,ZWSCCSS}.

The main purpose of this paper is to investigate the Segal-Bargmann transform $\mathcal{B}$ of Clifford algebra-valued functions, which has also been the focus of D. Pe{\~n}a Pe{\~n}a, I. Sabadini and F. Sommen \cite{PSS}. We will define and examine the Segal-Bargmann module $\mathcal{F}^2(\CC^m,\C_m^{\CC})$, a higher-dimensional analogue of the classical Segal-Bargmann space.


It is known that there is a close relationship between the Gabor transform (short-time Fourier transform with a Gaussian window) and the Segal-Bargmann transform. Recently, this connection has been used to filter a signal embedded in white noise \cite{Flandrin,AHKR}. Therefore, we investigate the mapping properties of the Segal-Bargmann transform in the context of Clifford estimators.

We prove that $\mathcal{B}$ is a unitary operator up to a scaling constant, and that it maps an orthonormal basis of $L^2(\RR^m,\C_m^{\RR})$ to an orthonormal basis of the Segal-Bargmann module $\mathcal{F}^2(\CC^m,\C_m^{\CC})$. For that, we will use Sommen's Clifford-Hermite functions $\left\{\phi_{l,k,j}\right\}$ as an $L^2$ basis.

We also lay out that the Segal-Bargmann transform can be expanded to a series $\big(\mathcal{B}f\big)(z)=\sum\limits_{l=0}^\infty\sum\limits_{k=0}^\infty\sum\limits_{j=1}^{\mathrm{dim}(M_l^+(k))}\Psi_{l,k,j}(\uz)\langle\phi_{l,k,j},f\rangle$ with a dictionary $\left\{\Psi_{l,k,j}\right\}$ of the Segal-Bargmann module and that this series converges absolutely locally uniformly.

The paper is organised as follows. In Section \ref{subsec:Cliffordalgebras} and \ref{subsec:HilbertCliffmod} we give an overview of basic Clifford analysis and of Hilbert Clifford-modules, which replace Hilbert spaces in our context. Section \ref{sub:Pk} deals with a certain class of Clifford-valued functions, the inner spherical monogenics, which are central to the construction of a basis for the function spaces that we deal with. In section \ref{sub:STFT}, we present the short-time Fourier transform as in important tool for our work.

After we have established these preliminary notes, we introduce Sommen's generalized Clifford Hermite polynomials and their relevant properties in section \ref{sec:CliffHemitepol}. In section \ref{sec:BargmannTransform}, we formally introduce the Segal-Bargmann transform and the Segal-Bargmann space of the classical, non-Clifford case, before we establish its analogue, the Segal-Bargmann module, in section \ref{sec:BargmannModules} and show some important properties of the Segal-Bargmann transform of Clifford algebra-valued functions. We conclude our paper with section \ref{sec:dictionary} by constructing a dictionary $\left\{\Psi_{l,k,j}\right\}$ for the Segal-Bargmann transform and  proving the convergence of the series representation $\sum\limits_{l=0}^{\infty}\sum\limits_{k=0}^\infty\sum\limits_{j=1}^{\mathrm{dim}(M_l^+(k))}\Psi_{l,k,j}(\uz)\langle\phi_{l,k,j},f\rangle$.

\section{Preliminaries}

\label{sec:prelim}

\subsection{Clifford algebras}

\label{subsec:Cliffordalgebras}

While real Clifford algebras have gained much interest in mathematical research since W. Clifford wrote about them in 1878, cf. \cite{Funktionentheorie}, complex Clifford algebras are a fairly recent topic of interest. In our work, we deal with both cases. We take notations and properties mainly from \cite{BDS-1982}, in which the real version is displayed, and adopt them to fit the complex case. For that, we work close to J. Ryan's \textit{Complexified clifford analysis} \cite{Ryan}, in which a detailed extension of real to complex Clifford algebras is developed.

We will write $\mathbb{N}=\left\{1,2,3,\dots\right\}$ and $\mathbb{N}_0=\left\{0,1,2,\dots\right\}$. Let $n\in\mathbb{N}_0$ and $\C_n^{\RR}$ denote the real Clifford algebra over $\RR^n$ and $\C_n^{\CC}$ the complex Clifford algebra over $\CC^n$. Both are based on the multiplication rules

$$ \begin{array}{cl} e_ie_j + e_je_i = 0, & i\not= j, \\ e_j^2 = -1, & i=1,2,\ldots , n. \end{array} $$

and have $e_0 \equiv 1$ as their unit element.

An arbitrary element of $\C_n^{\RR}$ or $\C_n^{\CC}$ is called a Clifford number and is given by

$$ a = \sum_{A} a_{A}e_{A},$$

where $a_{A}\in\RR$ or $a_{A}\in\CC$, resp., and for each $A=(n_1,\dots,n_l)$ with \linebreak $1\leq n_1 < n_2 < \ldots < n_l \leq m$, it is $e_{A} = e_{n_1}e_{n_2}\dots e_{n_l}$. The coefficient $a_0$ is called the scalar part of $a$ and $\underline{a} = \sum_{j=1}^n a_je_j$ a Clifford vector.

Similar to the complex conjugation $\overline{\phantom{x}}^{\CC}$, we can define involutions $\overline{\phantom{x}}$ for the real and $^\dagger$ for the complex Clifford algebra. Let 

\[\bar{e}_A=\left(-1\right)^{\frac{|A|(|A|+1)}{2}}e_A.\]

Then 

\[\bar{a}=\sum\limits_Aa_A\bar{e}_A\]

for $a\in\C_n^{\RR}$, and 

\[a^{\dagger}=\sum\limits_A\overline{a}_A^{\CC}\bar{e}_A\]\label{dagger}

for $a\in\C_n^{\CC}$.

We refer to \cite{GilbertMurray} and state that $\C_n^{\RR}$ becomes a finite dimensional Hilbert space with the inner product

$$ (a,b)_0 = [\overline{a}b]_0 = \sum_A a_Ab_A$$ 

for all $a,b \in \C_n^{\RR}$,

and has Hilbert space norm

$$ |a|_0 = \sqrt{ (a, a )_0} = \sqrt{\sum_A |a_A|^2}. $$

The inner product on $\C_n^{\RR}$ extends to a sesqui-linear inner product

$$ (a,b)_0 = [a^{\dagger}b]_0 = \sum_A \overline{a_A}^{\CC}b_A$$

for $a,b \in \C_n^{\CC}$.

It can be shown that Clifford algebras are $C^*$-algebras, see \cite{GilbertMurray}.

\begin{prop}

Under the involution $a \to a^{\dagger}$ each $\C_n^{\CC}$ is a

complex $C^*$-algebra which is a complexification of the real $C^*$-algebra $\C_n^{\RR}.$

\end{prop}

\subsection{Hilbert Clifford-modules}

\label{subsec:HilbertCliffmod}

We want to consider spaces of $\C_m^{\RR}$- or $C_m^{\CC}$-valued functions.  For that purpose, we need an anologue to the classical $L^2$ spaces. Since the elements of a Clifford algebra do not form a field, we work in Clifford-modules. The following two definitions are taken from \cite{BDS-1982} and adapted for the complex case; the real case is contained implicitly.

\begin{defn} $X_{(r)}$ is a \textbf{unitary right $\C_n^{\CC}$-module}, when $(X_{(r)}, +)$ is an abelian group and the mapping $(f, a) \to fa$ from $ X_{(r)}\times\C_n^{\CC} \to X_{(r)}$ is defined such that for all $a,b\in \C_n^{\CC}$ and $f,g\in X_{(r)}:$

\begin{enumerate}

\item $f(a + b) = fa + fb,$

\item $f(ab) = (fa)b, $

\item $(f+g)a = fa + ga, $

\item $fe_0 = f.$

\end{enumerate}

\end{defn}

We define an inner product on a unitary right $\C_n^{\CC}$-module as follows.

\begin{defn}

Let $H_{(r)}$ be a unitary right $\C_n^{\CC}$-module. Then a function\linebreak

$\langle\cdot, \cdot \rangle: H_{(r)}\times H_{(r)} \to \C_n^{\CC} $ is an \textbf{inner product} on $H_{(r)}$ if for all $f,g,h \in H_{(r)}$ and $a \in \C_n^{\CC},$ 

\begin{enumerate}

\item $\langle f, g+h\rangle = \langle f,g\rangle + \langle f,h\rangle$

\item $\langle f,ga\rangle = \langle f,g\rangle a $

\item $\langle f,g\rangle  = \langle g,f\rangle^{\dagger}$

\item $\langle f,f\rangle _0\in\RR_0^+$ and $\langle f,f\rangle _0=0$ if and only if $f=0$

\item $\langle fa,fa\rangle _0\leq|a|_0^2\langle f,f\rangle _0.$

\end{enumerate}

The accompanying \textbf{norm} on $H_{(r)}$ is $\left\|f\right\|^2=\langle f,f\rangle_0$.

\end{defn}

We now give an important property of the inner product.

\begin{prop}\label{prop:fg0}

If $\langle\cdot,\cdot\rangle$ is an inner product on a unitary right $\C_n^{\CC}$-module $H_{(r)}$ and $\left\|f\right\|^2=\langle f,f\rangle_0$ then

\[\left|\langle f,g\rangle\right|_0\leq 2^n\left\|f\right\|\,\left\|g\right\|\]

for all $f,g\in H_{(r)}$.

\end{prop}

Proof:

We use the definition of the norm on $\C_n^{\CC}$, $|a|_0^2=\sum\limits_A |a_A|^2$, and the fact that

\begin{equation} 
\left[ae_A\right]_0=\left[\sum\limits_Ba_Be_Be_A\right]_0=\left[a_Ae_Ae_A\right]_0=-a_A\label{eq:aeA}
\end{equation}

for all $a\in\C_n^{\CC}$. Also, if we consider $H_{(r)}$ to be a vector space over $\CC$ with inner product $\langle\cdot,\cdot\rangle_0$, we know that the Cauchy-Schwartz inequality

\begin{eqnarray}\label{eq:Norm-CS}
\left|\langle f,g\rangle_0\right|^2\leq\langle f,f\rangle_0\cdot\langle g,g\rangle_0=\left\|f\right\|^2\,\left\|g\right\|^2
\end{eqnarray}

has to be true. Now, we get

\begin{align*}
\left|\langle f,g\rangle\right|_0^2&=\sum\limits_A\left|\langle f,g\rangle_A\right|^2 \stackrel{(\ref{eq:aeA})}{=}\sum\limits_A\left|\left[\langle f,g\rangle e_A\right]_0\right|^2\\
&\stackrel{(ii)}{=}\sum\limits_A\left|\langle f,ge_A\rangle_0\right|^2 \stackrel{(\ref{eq:Norm-CS})}{\leq}\sum\limits_A\|f\|^2\,\|ge_A\|^2\\
&\stackrel{(v)}{=}\sum\limits_A\|f\|^2\,\|g\|^2\cdot |e_A|_0^2 =\sum\limits_A\|f\|^2\,\|g\|^2\\
&=2^n\|f\|^2\,\|g\|^2. \quad \square
\end{align*}

As an analogue to Hilbert (vector) spaces, we now define Hilbert modules.

\begin{defn}

Let $H_{(r)}$ be a unitary right $\C_n^{\CC}$-module provided with an inner product $(\cdot, \cdot).$ Then it is called a \textbf{right Hilbert $\C_n^{\CC}$-module} if it is complete for the norm topology derived from the inner product.

\end{defn}

Let $m\in\mathbb{N}=\{1,2,3,\dots\}$. We now consider the unitary right $\C_m^{\RR}$-module of functions from $\RR^m$ to $\C_m^{\RR}$. A function ${f:\Omega\subset\RR^m\to\C_m^{\RR}}$ maps the vector variable $\ux = \sum_{j=1}^m x_je_j$ to a Clifford number and can be written as

\[f(\ux)=\sum\limits_Ae_Af_A(\ux),\]

where $f_A:\RR^m\to\RR$ \cite{BDS-1982}. We define an inner product as follows.

\begin{defn} 

Let $h$ be a positive function on $\RR^m.$ Then the inner product $\langle\cdot,\cdot\rangle_{L^2(\RR^m,h,\C_m^{\RR})}$ is defined as

$$\langle f,g\rangle_{L^2(\RR^m,h,\C_m^{\RR})} = \int_{\RR^m} \overline{f(\ux)} g(\ux) h(\ux)   d\ux,$$

where $ d\ux$ stands for the Lebesgue measure on $\RR^m,$ and the associated norm is $||f||_{L^2(\RR^m,h,\C_m^{\RR})}^2 = \left[\langle f, f \rangle_{L^2(\RR^m,h,\C_m^{\RR})}\right]_0 .$ 

\end{defn}

The unitary right Clifford-module of measurable functions on $\RR^m$ for which $||f||_{L^2(\RR^m,h,\C_m^{\RR})} < \infty $ is a right Hilbert Clifford-module, which we denote by $L^2(\RR^m, h, \C_m^{\RR})$. In this paper, we will focus on the case where $h(\ux)=1$. Then the right Hilbert Clifford-module will simply be denoted by $L^2(\RR^m,\C_m^{\RR})$ and the inner product by $\langle\cdot,\cdot\rangle_{L^2(\RR^m,\C_m^{\RR})}$.

We also work on functions with values in a complex Clifford algebra, i.e. $f:\Omega\subset\CC^m\to\C_m^{\CC}$. For $\uz=\sum_{j=1}^m z_je_j$, with complex $z_j$, $j=1,\dots,m$, we have

\[f(\uz)=\sum\limits_Ae_Af_A(\uz)\]

with $f_A:\CC^m\to\CC$. Analogously to the real case, we can define the right Hilbert Clifford-module $L^2(\CC^m,h,\C_m^{\CC})$, where $h$ is a positive function over $\CC^m$. Here,

\[\langle f,g\rangle_{L^2(\CC^m,h,\C_m^{\CC})} = \int_{\CC^m} f^\dagger(\uz) g(\uz) h(\uz) d\ux\,d\uy\]

with $\uz=\ux+i\uy$, where $\dagger$ denotes the involution on $\C_m^{\CC}$, cf. page \pageref{dagger}. The associated norm is $||f||_{L^2(\CC^m,h,\C_m^{\CC})}^2=\left[\langle f,f\rangle_{L^2(\CC^m,h,\C_m^{\CC})}\right]_0$. Particularly important to our work will be those spaces $L^2(\CC^m,h,\C_m^{\CC})$ for which $h$ is defined as the Gaussian function $h(\uz)=\frac{e^{-|\uz|^2/2}}{\pi^m}$, cf. section \ref{sec:BargmannModules}.

\begin{prop}\label{prop:normintnorm}

\begin{enumerate}

\item Let $f\in L^2(\RR^m,h,\C_m^{\RR}).$ Then

\[||f||_{L^2(\RR^m,h,\C_m^{\RR})}^2=\int\limits_{\RR^m}\left|f(\ux)\right|_0^2h(\ux) d\ux.\]

\item Let $f\in L^2(\CC^m,h,\C_m^{\CC}).$ Then

\[||f||_{L^2(\CC^m,h,\C_m^{\CC})}^2=\int\limits_{\CC^m}\left|f(\uz)\right|_0^2h(\uz) d\ux\,d\uy.\]

\end{enumerate}

\end{prop}

Proof:

We only show (ii) since the real case is equivalent, 

\begin{align*}
||f||_{L^2(\CC^m,h,\C_m^{\CC})}^2&=\left[\langle f,f\rangle_{L^2(\CC^m,h,\C_m^{\CC})}\right]_0\\
&=\int\limits_{\CC^m}\left[f^{\dagger}(\uz)f(\uz)\right]_0h(\uz) d\ux\,d\uy	=\int\limits_{\CC^m}\left|f(\uz)\right|_0^2h(\uz) d\ux\,d\uy. \quad \square
\end{align*}

\subsection{Inner spherical monogenics}

\label{sub:Pk}

Since many of the following results are similar for functions of real and complex Clifford algebras, we will state them for the real case.

Of particular importance, when dealing with Clifford algebra-valued functions, is the Dirac operator 

\[D_{\ux} = \sum_{j=1}^m e_j\partial_{x_j}\quad (D_{\uz} = \sum\limits_{j=1}^me_j\partial_{z_j}).\]

Left nullsolutions of $D_{\ux}$ ($D_{\uz}$) are called (complex) left monogenic functions.

Let $m\in\mathbb{N}$ and $\mathfrak{P}^s$ be the space of scalar-valued polynomials in $\RR^m$ ($\CC^m$). Then a Clifford polynomial is an element of $\mathfrak{P}^s \otimes \C_m^{\RR}$ ($\mathfrak{P}^s \otimes \C_m^{\CC}$).

An important class of polynomials are the so called (complex) inner spherical monogenics. A left inner spherical monogenic of order $k$ is a left monogenic homogeneous Clifford polynomial $P_k$ of degree $k$. The set of all left inner spherical monogenics of order $k$ is denoted by $M_l^+(k)$ and has the dimension \cite{BDSS}

\[\mathrm{dim}(M_l^+(k))=\binom{m+k-2}{k},\]

(with $\mathrm{dim}(M_l^+(0))=1$ for all $m\in\mathbb{N}$).

We will deal with inner spherical monogenics over both $\RR^m$ and $\CC^m$. To differentiate, we will write $P_k(\ux):\RR^m\to\C_m^{\RR}$ and $P_k(\uz):\RR^m\to\C_m^{\CC}$.

\subsection{Short-time Fourier Transform} 

\label{sub:STFT}

An important tool in time-frequency analysis is the short-time Fourier Transform. It allows to analyse a given signal simultaneously in the time and in the frequency domain, because it calculates the Fourier Transform not over the whole signal, but small blocks of it.

Given a signal $f(\ut)$ and a window function $\varphi(\ut)$, the short-time Fourier Transform $(V_{\varphi}f)(\ut,\uom)$ is classically defined as

\[(V_{\varphi}f)(\ut,\uom) = \frac{1}{\sqrt{2\pi}^m}\int\limits_{\RR^m}f(\ux)\overline{\varphi(\ux-\ut)}^{\CC}e^{-i\uom\cdot\ux}d\ux.\]

A commonly used window function is the Gaussian window because it provides a very good resolution of the studied signal \cite{groechenig}. 

It is given by $h(\ux)=e^{-\frac{|\ux|^2}{4}}$.

\section{Clifford Hermite polynomials}

\label{sec:CliffHemitepol}

We will now consider Clifford Hermite polynomials as a special class of Clifford polynomials.

In the classical case, the Hermite polynomials over $\RR$ can be obtained from the Taylor expansion of the function $z\mapsto e^{z^2/2}$,

$$e^{z^2/2} = \sum_{n=0}^{\infty} e^{x^2/2} \frac{t^n}{n!}H_n(ix).$$

They can also be calculated explicitly by

\[H_n(x)=(-1)^ne^{x^2}\frac{ d^n}{ dx^n}e^{-x^2}.\]

Through a similar expansion for $\RR^m$, F. Sommen defined \textit{radial Hermite polynomials} \cite{So}, which are explicitly given by 

$$H_{k,m}(\ux) = (-1)^ke^{\frac{|\ux|^2}{2}}D_{\ux}^k e^{-\frac{|\ux|^2}{2}}.$$

Since the radial Hermite polynomials only form a basis for a certain kind of $L^2$ functions, i.e. such functions that are defined on the real line, Sommen developed a more complex set of polynomials, starting from the monogenic extension of $e^{-|\ux|^2/2}P_k(\ux)$, where $P_k(\ux)$ is a left inner spherical monogenic of degree $k$, cf. section \ref{sub:Pk}. This lead him to what he called the generalized Hermite polynomials, which can be used to construct a basis of $L^2(\RR^m,\C_m^{\RR})$.

\begin{defn}\label{def:CliffHerm} The \textbf{generalized Clifford Hermite polynomials} $H_{l,m,k}$, $l,k\in\mathbb{N}_0$, are given by

\begin{equation}\label{eq:CliffHermitePolynomials}
H_{l,m,k}P_k(\ux) = e^{-\frac{|\ux|^2}{2}}(-1)^lD_{\ux}^l\left(e^{-\frac{|\ux|^2}{2}}P_k(\ux)\right)
\end{equation}

where $P_k(\ux)$ is a left inner spherical monogenic of degree $k$.

\end{defn}

An important property of the generalized Clifford Hermite polynomials is their orthogonality \cite{BDSKS}.

\begin{thm}\label{th:CliffordHermiteOrtho}

Let $H_{l,m,k_1}$ and $H_{t,m,k_2}$ be generalized Clifford Hermite polynomials and $P_{k_1}(\ux)$ and $P_{k_2}(\ux)$ inner spherical monogenics of order $k_1$ and $k_2$, resp. Then

$$\int_{\RR^m} e^{-\frac{|\ux|^2}{2}}\, \overline{H_{l,m,k_1}(\ux)P_{k_1}(\ux)}\, H_{t,m,k_2}(\ux)P_{k_2}(\ux)\,  d\ux = \gamma_{l,k_1} \delta_{l,t} \delta_{k_1,k_2},$$

with

\begin{align*}
\gamma_{2p,k}&=\frac{2^{2p+m/2+k}p!\sqrt{\pi}^m\Gamma\left(\frac{m}{2}+k+p\right)}{\Gamma\left(\frac{m}{2}\right)},\\
\gamma_{2p+1,k}&=\frac{2^{2p+m/2+k+1}p!\sqrt{\pi}^m\Gamma\left(\frac{m}{2}+k+p+1\right)}{\Gamma\left(\frac{m}{2}\right)}.\\
\end{align*}

\end{thm}

Building on the orthogonality, Sommen and his colleagues established an orthonormal basis of $L^2(\RR^m,\C_m^{\RR})$, cf. \cite{So,BDSKS}.

\begin{thm}\label{th:L2basis}

Let $\gamma_{l,k}$, $k,l\in\mathbb{N}_0$, be as defined in Theorem \ref{th:CliffordHermiteOrtho}. For each $k\in\mathbb{N}_0$, let further  $\left\{P_k^{(j)}(\ux)\right\}_{j=1,2,\dots,\mathrm{dim}(M_l^+(k))}$ be an orthonormal basis of $M_l^+(k)$. 

\begin{equation}\label{eq:L2basis}
\left\{\frac{1}{\sqrt{\gamma_{l,k}}}H_{l,m,k}(\ux)P_k^{(j)}(\ux)e^{-\frac{|\ux|^2}{4}}:l,k\in\mathbb{N}_0,j\leq\mathrm{dim}(M_l^+(k))\right\}
\end{equation}

forms an orthonormal basis of $L^2(\RR^m, \C_m^{\RR})$. 

\end{thm}

Each element of (\ref{eq:L2basis}) depends on $l$, $k$ and the chosen basis of $M_l^+(k)$, which contains $\mathrm{dim}(M_l^+(k))=\binom{m+k-2}{k}$ elements, cf. section \ref{sub:Pk}.

\section{Segal-Bargmann transform}

\label{sec:BargmannTransform}

The first very general version of the Segal-Bargmann space goes back to V. Fock's theory of 1932 of the quantum state space of particles \cite{Fock}. Here, we will consider the more specific finite-dimensional version of the following definitions taken from \cite{PSS}. We will transfer those definitions to the Clifford case in section \ref{sec:BargmannModules}.

\begin{defn} 

The \textbf{Segal-Bargmann space} $\mathcal{F}^2(\CC^m,\CC)$ is defined as the \linebreak Hilbert space of entire functions $f$ in $\CC^m$ which are square-integrable with respect to the $2m$-dimensional Gaussian density, i.e.,

$$ \frac{1}{\pi^m} \int_{\CC^m} e^{-|\uz|^2} |f(\uz)|^2 \,  d\ux\,d\uy < \infty, \quad \uz = \ux + i\uy. $$

It is equipped with the inner product

$$ \langle f,g \rangle_{\mathcal{F}^2(\CC^m,\CC)} =  \frac{1}{\pi^m} \int_{\CC^m} e^{-|\uz|^2} \overline{f(\uz)} g(\uz) \,  d\ux\,d\uy. $$

\end{defn}

The Segal-Bargmann transform connects the Bargmann space with the Hilbert space $L^2(\RR^m,\RR)$ by mapping the ladder onto the former.

\begin{defn}

The \textbf{Segal-Bargmann transform} $\mathcal{B}$ from $L^2(\RR^m,\RR)$ to\linebreak $\mathcal{F}^2(\CC^m,\CC)$ is defined by

\begin{equation}\label{def:BargmannTransform}
\big(\mathcal{B}f\big)(\uz) = \frac{1}{\sqrt{2\pi}^m} \int_{\RR^m} e^{-\frac{\uz\cdot\uz}{2} + \ux \cdot\uz - \frac{\ux\cdot\ux}{4}} f(\ux)\, d\ux,
\end{equation}

with $\ux \cdot \uz = \sum_{j=1}^m x_jz_j,$ for any $f\in L^2(\RR^m,\RR).$

\end{defn}

The Segal-Bargmann transform is a linear operator. It can also be expressed in terms of a short-time Fourier Transform (cf. section \ref{sub:STFT}).

\begin{prop}\label{prop:STFT-B}

Let $\mathcal{B}$ be the Segal-Bargmann transform and $V_\varphi$ the short-time Fourier Transform with window $\vphi(\ux)=e^{-\frac{|\ux|^2}{4}}.$ Then for all $f\in L^2(\RR^m,\RR)$,

\[(V_{\vphi}f)(2\ut,-\uom)= e^{-\frac{|\uz|^2}{2}}e^{i\ut\cdot \uom}\big(\mathcal{B}f\big)(\uz),\quad \uz = \ut +i\uom.\]

\end{prop}

Proof:  

\begin{align*}
(V_\varphi f)(2\ut,-\uom) & = \frac{1}{\sqrt{2\pi}^m} \int_{\RR^m} f(\ux) e^{-\frac{|\ux-2\ut|^2}{4}} e^{i\uom\cdot\ux}  d\ux \\
 & = \frac{1}{\sqrt{2\pi}^m} \int_{\RR^m} f(\ux)e^{-\frac{|\ux|^2}{4}+\ux\cdot\ut -|\ut|^2}e^{i\uom\cdot \ux}  d\ux  \\
& = \frac{1}{\sqrt{2\pi}^m} e^{-\frac{|\ut|^2}{2}}e^{-\frac{|\uom|^2}{2}}e^{i\ut\cdot\uom}\int_{\RR^m} f(\ux) e^{-\frac{|\ux|^2}{4} +\ux\cdot(\ut+i\uom) - \frac{(\ut+i\uom)^2}{2}}  d\ux \\
& = \frac{1}{\sqrt{2\pi}^m} e^{-\frac{|\ut|^2}{2}}e^{-\frac{|\uom|^2}{2}}e^{i\ut\cdot\uom}\int_{\RR^m} f(\ux) e^{-\frac{|\ux|^2}{4} + \ux\cdot \uz -\frac{\uz\cdot\uz}{2}}  d\ux \\
& =  e^{-\frac{|\uz|^2}{2}}e^{i\ut\cdot\uom}\big(\mathcal{B}f\big)(\uz) \quad \square
\end{align*}

A well-known property of the Segal-Bargmann transform is that it is a unitary operator up to a scaling constant.

\begin{prop}\label{prop:isometry-classical}

Let $\mathcal{B}:L^2(\RR^m,\RR)\to\mathcal{F}^2(\CC^m,\CC)$ be the Segal-Bargmann transform (\ref{def:BargmannTransform}). Then

\[\langle\mathcal{B}f,\mathcal{B}g\rangle_{\mathcal{F}^2(\CC^m,\CC)}=\frac{1}{\sqrt{2\pi}^m}\langle f,g\rangle_{L^2(\RR^m,\RR)}.\]

\end{prop}

Proof:

We use Proposition \ref{prop:STFT-B}, i.e.

\[\big(\mathcal{B}f\big)(\uz)= e^{\frac{|\uz|^2}{2}}e^{-i\ut\cdot \uom}(V_{\vphi}f)(2\ut,-\uom),\]

with $\uz=\ut+i\uom$ and $\vphi(\ux)=e^{-\frac{|\ux|^2}{4}}.$ Then

\begin{align*}
\langle\mathcal{B}f,\mathcal{B}g\rangle_{\mathcal{F}^2(\CC^m,\CC)}
&=\frac{1}{\pi^m}\int\limits_{\CC^m}\overline{\big(\mathcal{B}f\big)(\uz)}^{\CC}\big(\mathcal{B}g\big)(\uz)e^{-|\uz|^2} d\ux\,d\uy\\
&=\frac{1}{\pi^m}\int\limits_{\CC^m}e^{\frac{|\uz|^2}{2}}e^{i\ut\cdot\uom}\overline{\big(V_\varphi f\big)(2\ut,-\uom)}^{\CC}\\
&\quad\quad\quad\quad\cdot e^{\frac{|\uz|^2}{2}}e^{-i\ut\cdot\uom}\big(V_\varphi g\big)(2\ut,-\uom)e^{-|\uz|^2} d\uom\, d\ut\\
&=\frac{1}{\pi^m}\int\limits_{\CC^m}\overline{\big(V_\varphi f\big)(2\ut,-\uom)}^{\CC}\big(V_\varphi g\big)(2\ut,-\uom) d\uom\, d\ut\\
&=\frac{1}{\pi^m}\int\limits_{\CC^m}\overline{\frac{1}{\sqrt{2\pi}^m}\int\limits_{\RR^m}f(\ux)e^{-\frac{|\ux-2\ut|^2}{4}}e^{i\uom\cdot\ux} d\ux}^{\CC}\\
&\quad\quad\quad\quad\cdot\frac{1}{\sqrt{2\pi}^m}\int\limits_{\RR^m}g(\ux)e^{-\frac{|\ux-2\ut|^2}{4}}e^{i\uom\cdot\ux} d\ux\, d\uom\, d\ut.
\end{align*}

Let $\varphi(\cdot-2\ut)$ denote the Gaussian window translated by $-2\ut$ and $\mathcal{F}$ the Fourier Transform in $\RR^m$. Thus

\[\langle\mathcal{B}f,\mathcal{B}g\rangle_{\mathcal{F}^2(\CC^m,\CC)} =\frac{1}{\pi^m}\int\limits_{\CC^m}\overline{\mathcal{F}^{-1}\big(f\cdot\varphi(\cdot-2\ut)\big)(\uom)}^{\CC}\mathcal{F}^{-1}\big(g\cdot\varphi(\cdot-2\ut)\big)(\uom) d\uom\; d\ut.\]

The Plancherel Theorem now gives us

\begin{align*}
\langle\mathcal{B}f,\mathcal{B}g\rangle_{\mathcal{F}^2(\CC^m,\CC)} &=\frac{1}{\pi^m}\int\limits_{\RR^m}\int\limits_{\RR^m}\overline{f(\ux)e^{-\frac{|\ux-2\ut|^2}{4}}}^{\CC}g(\ux)e^{-\frac{|\ux-2\ut|^2}{4}} d\ut\, d\ux\\
&=\frac{1}{\pi^m}\int\limits_{\RR^m}f(\ux)g(\ux)\int\limits_{\RR^m}\left(e^{-\frac{|\ux-2\ut|^2}{4}}\right)^2 d\ut\, d\ux.
\end{align*}

Last, we substitute $\underline{u}=2\ut-\ux$, and use the fact that $\int\limits_{\RR^m}e^{-\frac{|\underline{u}|^2}{2}} d\underline{u}=\sqrt{2\pi}^m$. So,

\begin{align*}
\langle\mathcal{B}f,\mathcal{B}g\rangle_{\mathcal{F}^2(\CC^m,\CC)} &=\frac{1}{\pi^m}\int\limits_{\RR^m}f(\ux)g(\ux) d\ux\int\limits_{\RR^m}e^{-\frac{|u|^2}{2}}\frac{ d\underline{u}}{2^m}\\
&=\frac{1}{\sqrt{2\pi}^m}\langle f,g\rangle_{L^2(\RR^m,\RR)}. \quad \square
\end{align*}

Another important property of the Segal-Bargmann transform is its invertibility.

\vspace*{1cm}

\begin{prop}\label{prop:Binvertible}

\leavevmode

\begin{enumerate}

\item $F^2(\CC^m,\CC)$ is the image of $L^2(\RR^m,\RR)$ under the Segal-Bargmann transform.

\item The Segal-Bargmann transform is invertible.

\end{enumerate}

\end{prop}

Proof:

For the proof of (i) we refer to \cite{groechenig}. (ii) then follows directly from the fact that the transform is unitary up to a constant, cf. Proposition \ref{prop:isometry-classical}. $\square$

\section{Segal-Bargmann modules}

\label{sec:BargmannModules}

We will now look at how the Segal-Bargmann transform of Definition \ref{def:BargmannTransform} acts on Clifford algebra-valued functions. So, from now on, let $f$ be an element of $L^2(\RR^m,\C_m^{\RR})$. Then,

\[\big(\mathcal{B}f\big)(\uz) = \frac{1}{\sqrt{2\pi}^m} \int_{\RR^m} e^{-\frac{\uz\cdot\uz}{2} + \ux \cdot\uz - \frac{\ux\cdot\ux}{4}} f(\ux)\, d\ux,\]

is a function with values in the complex Clifford algebra, $\mathcal{B}f:\CC^m\to\C_m^{\CC}$. Note that Proposition \ref{prop:STFT-B} holds for functions of $L^2(\RR^m,\C_m^\RR)$ as well.

Consider the function space 

$L^2(\CC^m,\frac{e^{-|\uz|^2}}{\pi^m},\C_m^{\CC})$ 


as defined in section \ref{subsec:HilbertCliffmod}. Just as in the real case, the Segal-Bargmann transform of Clifford algebra-valued functions is unitary up to a scaling constant, as the following proposition shows.

\begin{prop}\label{prop:SegalBargmannIsometry} If $f\in L^2(\RR^m,\C_m^{\RR}),$ then 

\[\langle \mathcal{B}f, \mathcal{B}g \rangle_{L^2(\CC^m,\frac{e^{-|\uz|^2}}{\pi^m},\C_m^{\CC})}=\frac{1}{\sqrt{2\pi}^m}\langle f, g \rangle_{L^2(\RR^m, \C_m^{\RR})}.\]

\end{prop}

Proof:

Since the Segal-Bargmann transform is linear, $f=\sum\limits_A f_Ae_A$ implies that $\mathcal{B}f=\sum\limits_A \mathcal{B}f_Ae_A$. Hence,

\begin{align*}
\langle \mathcal{B}f, \mathcal{B}g \rangle_{L^2(\CC^m,\frac{e^{-|\uz|^2}}{\pi^m},\C_m^{\CC})} & = \frac{1}{\pi^m}\int_{\CC^m} \big(\mathcal{B}f\big)^{\dagger}(\uz) \big(\mathcal{B}g\big)(\uz) e^{-|\uz|^2} d\ux\,d\uy   \\
 & = \sum_{A,B} \frac{1}{\pi^m}\int\limits_{\CC^m}\overline{\big(\mathcal{B}f_A\big)}^{\CC} (\uz)\overline{e_A} \big(\mathcal{B}g_B\big)(\uz)e_B e^{-|\uz|^2} d\ux\,d\uy  \\
 & = \sum_{A,B} \langle \mathcal{B}f_A, \mathcal{B}g_B \rangle_{\mathcal{F}^2(\CC^m,\CC)} \overline{e_A}e_B
\end{align*}

In Proposition \ref{prop:isometry-classical} we have shown that 

$\langle \mathcal{B}f, \mathcal{B}g \rangle_{\mathcal{F}^2(\CC^m, \CC)}=\frac{1}{\sqrt{2\pi}^m}\langle f, g \rangle_{L^2(\RR^m, \RR)}$

 is true for the classical Segal-Bargmann transform $\mathcal{B}: L^2(\RR^m,\RR)\to\mathcal{F}^2(\CC^m,\CC)$. Therefore

\begin{align*}
\langle \mathcal{B}f, \mathcal{B}g \rangle_{L^2(\CC^m,\frac{e^{-|\uz|^2}}{\pi^m},\C_m^{\CC})}  & = \sum_{A,B} \frac{1}{\sqrt{2\pi}^m}\langle f_A, g_B \rangle_{L^2(\RR^m,\RR)} \overline{e_A}e_B \\
 & = \frac{1}{\sqrt{2\pi}^m}\sum\limits_{A,B}\,\int\limits_{\RR^m}\overline{f_A(\ux)}g_B(\ux)\overline{e_A}e_B d\ux\\
 & = \frac{1}{\sqrt{2\pi}^m}\int\limits_{\RR^m}\overline{f(\ux)}g(\ux) d\ux  = \frac{1}{\sqrt{2\pi}^m}\langle f, g \rangle_{L^2(\RR^m, \C_m^{\RR})}
 \quad \square
\end{align*}

A direct consequence is the following corollary.

\begin{cor}\label{cor:SegalBargmannIsometry}

Let $\left\|\cdot\right\|_{L^2(\CC^m,\frac{e^{-|\uz|^2}}{\pi^m},\C_m^{\CC})}=\sqrt{\left[\langle\cdot,\cdot\rangle_{L^2(\CC^m,\frac{e^{-|\uz|^2}}{\pi^m},\C_m^{\CC})}\right]_0}$. Then

\[\left\|\mathcal{B}f\right\|_{L^2(\CC^m,\frac{e^{-|\uz|^2}}{\pi^m},\C_m^{\CC})}^2=\frac{1}{\sqrt{2\pi}^m}\left\|f\right\|_{L^2(\RR^m,\C_m^{\RR})}^2.\]

Thus $\mathcal{B}$ is an isometry from $L^2(\RR^m,\C_m^{\RR})$ into $L^2(\CC^m,\frac{e^{-|\uz|^2}}{\pi^m},\C_m^{\CC})$ up to 

$\frac{1}{\sqrt{2\pi}^m}$.

\end{cor}







In Theorem \ref{th:zlPk}, an orthonormal basis of the space $L^2(\RR^m,\C_m^{\RR})$ was established. The following theorem shows that the Segal-Bargmann transform maps the elements of this basis onto functions $\uz^lP_k(\uz)$.

\begin{thm}\label{th:zlPk}

Let $\mathcal{B}$ be the Segal-Bargmann transform, $H_{l,m,k}$ a generalized Clifford Hermite Polynomial as defined in Definition \ref{def:CliffHerm} and $P_k$ an inner spherical monogenic of degree $k$. Then 

\[\left(\mathcal{B}\left(H_{l,m,k}(\ux)e^{-\frac{|\ux|^2}{4}}P_k(\ux)\right)\right)(\uz)=\uz^lP_k(\uz)\]

\end{thm}

Proof:

Our first step follows \cite{PSS}. Here, 

\begin{align*}
&\left(\mathcal{B}\left(H_{l,m,k}(\ux)e^{-\frac{|\ux|^2}{4}}P_k(\ux)\right)\right)(\uz)\\
&\quad\quad\quad\quad=\frac{1}{\sqrt{2\pi}^m}\int\limits_{\RR^m}e^{-\frac{\uz\cdot\uz}{2}+\ux\cdot\uz-\frac{\ux\cdot\ux}{2}}H_{l,m,k}(\ux)P_k(\ux) d\ux\\
&\quad\quad\quad\quad\stackrel{(\ref{eq:CliffHermitePolynomials})}{=}\frac{(-1)^l}{\sqrt{2\pi}^m}\int\limits_{\RR^m}e^{-\frac{\uz\cdot\uz}{2}+\ux\cdot\uz}D_{\ux}^l\left(e^{-\frac{|\ux|^2}{2}}P_k(\ux)\right) d\ux\\
&\quad\quad\quad\quad=\frac{1}{\sqrt{2\pi}^m}\int\limits_{\RR^m}D_{\ux}^l\left(e^{-\frac{\uz\cdot\uz}{2}+\ux\cdot\uz}\right)e^{-\frac{|\ux|^2}{2}}P_k(\ux)d\ux\\
&\quad\quad\quad\quad=\frac{1}{\sqrt{2\pi}^m}\uz^l\int\limits_{\RR^m}e^{-\frac{\uz\cdot\uz}{2}+\ux\cdot\uz-\frac{\ux\cdot\ux}{4}}P_k(\ux)e^{-\frac{|\ux|^2}{4}}d\ux\\
&\quad\quad\quad\quad=\uz^l\left(\mathcal{B}\left(P_k(\ux)e^{-\frac{|\ux|^2}{4}}\right)\right)(\uz).
\end{align*}

Next, we calculate $\mathcal{B}(P_k(\ux)e^{-|\ux|^2/4}) $ using the windowed Fourier transform. We obtain

\begin{align*}
V_{\vphi}\left(P_k(\ux)e^{-\frac{|\ux|^2}{4}}\right)(2\ut, -\uom) & = \frac{1}{\sqrt{2\pi}^m} \int_{\RR^m} P_k(\ux)e^{-\frac{|\ux|^2}{4}}e^{-\frac{|\ux-2\ut|^2}{4}}e^{i\uom\cdot\ux} d\ux \\
 & =  \frac{1}{\sqrt{2\pi}^m} P_k\left(-i\partial_{\uom}\right) \int_{\RR^m} e^{-\frac{|\ux|^2}{4}}e^{-\frac{|\ux|}{4}+\ux\cdot\ut -|\ut|^2}e^{i\uom\cdot\ux} d\ux \\
 & =  \frac{1}{\sqrt{2\pi}^m} P_k\left(-i\partial_{\uom}\right) \int_{\RR^m} e^{-\frac{|\ux|}{2}+\ux\cdot\ut -\frac{|\ut|^2}{2}}e^{-\frac{|\ut|}{2}}e^{i\uom\cdot\ux} d\ux \\
 & = \frac{1}{\sqrt{2\pi}^m} e^{-\frac{|\ut|^2}{2}}P_k\left(-i\partial_{\uom}\right) \int_{\RR^m} e^{-\frac{|\ux-\ut|^2 }{2}}e^{i\uom\cdot\ux} d\ux \\
 &  = \frac{1}{\sqrt{2\pi}^m} e^{-\frac{|\ut|^2}{2}}P_k\left(-i\partial_{\uom}\right)\left(e^{i\uom\cdot\ut} \int_{\RR^m} e^{-\frac{|\ux|^2 }{2}}e^{i\uom\cdot\ux} d\ux \right)
\end{align*}

Since $\frac{1}{\sqrt{2\pi}^m}\int_{\RR^m}e^{-\frac{|\ux|^2}{2}}e^{i\uom\cdot\ux}\, d\ux$ is the inverse Fourier Tranform of $e^{-\frac{|\uom|^2}{2}}$, which is an invariant, we get

\begin{align*}
V_{\vphi}\left(P_k(\ux)e^{-\frac{|\ux|^2}{4}}\right)(2\ut, -\uom) & = e^{-\frac{|\ut|^2}{2}}P_k\left(-i\partial_{\uom}\right)\left(e^{i\uom\cdot\ut -\frac{|\uom|^2}{2}}\right) \\
 & = e^{-\frac{|\ut|^2}{2}}P_k\left(-i(i\ut-\omega)\right) \left(e^{i\uom\cdot\ut -\frac{|\uom|^2}{2}}\right) \\
 & = e^{-\frac{|\ut|^2}{2}} \left(e^{i\uom\cdot\ut -\frac{| \uom|^2}{2}}\right) P_k(\uz) \\
 & = e^{-\frac{|\uz|^2}{2}} e^{i\uom\cdot \ut}P_k(\uz) 
\end{align*}

with $\uz=\ut+i\uom$. Because of Proposition \ref{prop:STFT-B}, this leads to \[\left(\mathcal{B}\left(P_k(\ux)e^{-\frac{|\ux|^2}{4}}\right)\right)(\uz)=P_k(\uz).\]

Together with the first step, the proof is complete. $\square$

We can now define an analogue to the classical Segal-Bargmann space.

\begin{defn}

The closure of

\[\mathrm{span}\left\{\uz^lP_k^{(j)}(\uz)\middle|l,k\in\mathbb{N}_0,j=1,\dots,\mathrm{dim}(M_l^+(k))\right\}\]

is called \textbf{Segal-Bargmann module} $\mathcal{F}^2(\CC^m,\C_m^{\CC})$.

\end{defn}

\begin{rem} In this definition and what follows we drop the property that a function of the Segal-Bargmann module (or space) has to be an entire functions. That means we consider the Segal-Bargmann module just as a weighted $L^2$-module.
\end{rem}






A consequence of Theorem \ref{th:zlPk} is the following.

\begin{cor}\label{cor:F2basis}

For all $l,k\in\mathbb{N}_0$, let $\left\{P_k^{(j)}(\ux)\right\}_{j=1,2,\dots,\mathrm{dim}(M_l^+(k))}$ be an orthonormal basis of $M_l^+(k)$ and $\gamma_{l,k}$ defined as in Theorem \ref{th:CliffordHermiteOrtho}. Then  \[\left\{\sqrt{\frac{(2\pi)^m}{\gamma_{l,k}}}\uz^lP_k^{(j)}(\uz)\middle|l,k\in\mathbb{N}_0,j=1,\dots,\mathrm{dim}(M_l^+(k))\right\}\] is an orthonormal basis of the Segal-Bargmann module $\mathcal{F}^2(\CC^m,\C_m^{\CC})$.

\end{cor}

\begin{prop}

Since the Segal-Bargmann transform is linear, Theorem \ref{th:zlPk} shows that it maps an element

$\phi_{l,k,j}(\ux)=\frac{1}{\sqrt{\gamma_{l,k}}}H_{l,m,k}(\ux)P_k^{(j)}(\ux)e^{-\frac{|\ux|^2}{4}}$

of the orthonormal basis of $L^2(\RR^m,\C_m^{\RR})$ (see Theorem \ref{th:L2basis}) onto

\[\left(\mathcal{B}\phi_{l,k,j}\right)(\uz)=\frac{1}{\sqrt{\gamma_{l,k}}}\uz^lP_k^{(j)}(\uz).\]

The statement now follows directly from Proposition \ref{prop:SegalBargmannIsometry} and Corollary \ref{cor:SegalBargmannIsometry}, which say that 

$\left\|\mathcal{B}\phi_{l,k,j}\right\|_{L^2(\CC^m,\frac{e^{-|\uz|^2}}{\pi^m},\C_m^{\CC})}=\frac{1}{\sqrt{2\pi}^m}\left\|\phi_{l,k,j}\right\|_{L^2(\RR^m,\C_m^{\RR})}=\frac{1}{\sqrt{2\pi}^m}$

and $\mathcal{B}$ is unitary up to the scaling constant.$\square$

\begin{thm}

The Segal-Bargmann module is the image of $L^2(\RR^m,\C_m^{\RR})$ under the Segal-Bargmann transform, i.e. 

\[\mathcal{F}^2(\CC^m,\C_m^{\CC})=L^2\left(\CC^m,\frac{e^{-|\uz|^2}}{\pi^m},\C_m^{\CC}\right).\]

\end{thm}

Proof:

First, let $F\in\mathcal{F}^2(\CC^m,\C_m^{\CC})$. By construction there has to exist a function $f\in L^2(\RR^m,\C_m^{\RR})$ so that $\mathcal{B}f=F$. Since $\mathcal{B}$ is unitary up to a constant, we know that $F\in L^2\left(\CC^m,\frac{e^{-|\uz|^2}}{\pi^m},\C_m^{\CC}\right)$. Hence

$\mathcal{F}^2(\CC^m,\C_m^{\CC})\subseteq L^2\left(\CC^m,\frac{e^{-|\uz|^2}}{\pi^m},\C_m^{\CC}\right)$.

We now show the opposite inclusion. Let $F\in L^2\left(\CC^m,\frac{e^{-|\uz|^2}}{\pi^m},\C_m^{\CC}\right)$. Then $F$ can be written as $F=\sum_AF_Ae_A$ with $F_A:\C^m\to\CC$ for all $A$. Since

\begin{align*}
\left\|F\right\|_{L^2\left(\CC^m,\frac{e^{-|\uz|^2}}{\pi^m},\C_m^{\CC}\right)}&=\langle\sum_AF_Ae_A,\sum_BF_Be_B\rangle_0\\
&=\sum_A\left\|F_A\right\|_{L^2\left(\CC^m,\frac{e^{-|\uz|^2}}{\pi^m},\C_m^{\CC}\right)}=\sum_A\left\|F_A\right\|_{L^2\left(\CC^m,\frac{e^{-|\uz|^2}}{\pi^m},\CC\right)}
\end{align*}

is finite if and only if $\left\|F_A\right\|_{L^2\left(\CC^m,\frac{e^{-|\uz|^2}}{\pi^m},\CC\right)}$ is finite for every $A$, we know that $F_A\in L^2(\CC^m,\frac{e^{-|\uz|^2}}{\pi^m},\CC)=\mathcal{F}^2(\CC_m,\CC)$, cf. Proposition \ref{prop:Binvertible} (ii).

Proposition \ref{prop:Binvertible} (i) tells us that $\mathcal{B}:L^2(\RR^m,\RR)\to\mathcal{F}^2(\CC_m,\CC)$ is invertible, so for each $A$ there exists $f_A\in L^2(\RR^m,\RR)$ so that $\mathcal{B}f_A=F_A$. Since $\mathcal{B}$ is linear,

\[F=\sum_AF_Ae_A=\sum_A(\mathcal{B}f_A)e_A=\mathcal{B}\left(\sum_Af_Ae_A\right),\]

so there exists a function $\sum_Af_Ae_A=f\in L^2(\RR^m,\C_m^{\RR})$ such that $\mathcal{B}f=F$. Therefore $F\in\mathcal{F}^2(\CC^m,\C_m^{\CC})$. $\square$

\section{A dictionary for the Segal-Bargmann transform}

\label{sec:dictionary}

In this section, we want to give a series representation for the Segal-Bargmann transform $\mathcal{B}$ on the right Clifford-module $L^2(\RR^m,\C_m^{\RR})$. By demonstrating that this representation converges absolutely locally uniformly, we will show that $\mathcal{B}f$ is well-defined and can be represented in kernel form. We work close to R. Bardenet and A. Hardy \cite{Bardenet-Hardy}, who have have shown similar characteristics of the classical Segal-Bargmann transform on $L^2(\RR^m,\RR)$ and other transforms.

For the rest of this section, we will shorten our notation by writing\linebreak $L^2=L^2(\RR^m,\C_m^{\RR})$, $\mathcal{F}^2=\mathcal{F}^2(\CC^m,\C_m^{\CC})$, $\langle\cdot,\cdot\rangle_{\mathcal{F}^2}=\langle\cdot,\cdot\rangle_{L^2(\CC^m,\frac{e^{-|\uz|^2}}{\pi^m},\C_m^{\CC})}$ and $\|\cdot\|_{\mathcal{F}^2}=\|\cdot\|_{L^2(\CC^m,\frac{e^{-|\uz|^2}}{\pi^m},\C_m^{\CC})}$.

Since the set $\left\{\phi_{l,k,j}\right\}_{l,k\in\mathbb{N}_0,j\in\{1,\dots,\mathrm{dim}(M_l^+(k))\}}$ of  Hermite functions 

\begin{equation}
\phi_{l,k,j}(\ux)=\frac{1}{\sqrt{\gamma_{l,k}}}H_{l,m,k}(\ux)P_k^{(j)}(\ux)e^{-\frac{|\ux|^2}{4}}\label{eq:phi2}
\end{equation}

is a basis of $L^2$, see section \ref{sec:CliffHemitepol}, each Clifford algebra-valued square integrable function $f(\ux)$ can be expanded as

\[f(\ux)=\sum\limits_{l=0}^\infty\sum\limits_{k=0}^\infty\sum\limits_{j=1}^{\mathrm{dim}(M_l^+(k))}\phi_{l,k,j}(\ux)\langle\phi_{l,k,j},f\rangle_{L^2}.\]

Hence,

\begin{align}
&\big(\mathcal{B}f\big)(\uz)\nonumber\\
&\quad=\frac{1}{\sqrt{2\pi}^m}\int\limits_{\RR^m} \left(\sum\limits_{l=0}^\infty\sum\limits_{k=0}^\infty\sum\limits_{j=1}^{\mathrm{dim}(M_l^+(k))}\phi_{l,k,j}(\ux)\langle\phi_{l,k,j},f\rangle_{L^2}\right)e^{-\frac{\uz\cdot\uz}{2}+\ux\cdot\uz-\frac{\ux\cdot\ux}{4}}d\ux\nonumber\\
&\quad=\sum\limits_{l=0}^\infty\sum\limits_{k=0}^\infty\sum\limits_{j=1}^{\mathrm{dim}(M_l^+(k))}\left(\frac{1}{\sqrt{2\pi}^m}\int\limits_{\RR^m}\phi_{l,k,j}(\ux)e^{-\frac{\uz\cdot\uz}{2}+\ux\cdot\uz-\frac{\ux\cdot\ux}{4}}d\ux\right)\langle\phi_{l,k,j},f\rangle_{L^2}\nonumber\\
&\quad=\sum\limits_{l=0}^\infty\sum\limits_{k=0}^\infty\sum\limits_{j=1}^{\mathrm{dim}(M_l^+(k))}\big(\mathcal{B}\phi_{l,k,j}\big)(\uz)\langle\phi_{l,k,j},f\rangle_{L^2}\nonumber\\
&\quad=\sum\limits_{l=0}^\infty\sum\limits_{k=0}^\infty\sum\limits_{j=1}^{\mathrm{dim}(M_l^+(k))}\frac{1}{\sqrt{\gamma_{l,k}}}\uz^lP_k^{(j)}(\uz)\langle\phi_{l,k,j},f\rangle_{L^2}\nonumber\\
&\quad=\sum\limits_{l=0}^\infty\sum\limits_{k=0}^\infty\sum\limits_{j=1}^{\mathrm{dim}(M_l^+(k))}\Psi_{l,k,j}(\uz)\langle\phi_{l,k,j},f\rangle_{L^2}\label{eq:Psi}
\end{align}

with $\Psi_{l,k,j}(\uz)=\frac{1}{\sqrt{\gamma_{l,k}}}\uz^lP_k^{(j)}(\uz)$.

To be able to show convergence of the series expansion (\ref{eq:Psi}), we need the following two lemmas.

\begin{lem}\label{lem:Pk-z}

Let $P_s(\uz)=\sum\limits_{|\alpha|=s}a_\alpha\uz^\alpha$ be a homogeneous $\C_m^{\CC}$-polyomial of degree $s$, with $a_\alpha\in\C_m^{\CC}$ for all $|\alpha|=s$. Then

\begin{enumerate}

\item $\left\|P_s(\uz)\right\|^2_{\mathcal{F}^2}=\sum\limits_{|\alpha|=s}\left|a_{\alpha}\right|_0^2\alpha!$,

\item $\left|P_s(\uz)\right|_0^2\leq\frac{1}{s!}\left\|P_s(\uz)\right\|^2_{\mathcal{F}^2}\left|\uz\right|_0^{2s}$.

\end{enumerate}

\end{lem}

Proof: \begin{enumerate}

\item We have

\begin{align*}
\left\|P_s(\uz)\right\|_{\mathcal{F}^2}^2 &=\left[\langle P_s(\uz),P_s(\uz)\rangle_{\mathcal{F}^2}\right]_0 =\frac{1}{\pi^m}\int\limits_{\CC^m}\left[P_s^{\dagger}(\uz)P_s(\uz)\right]_0e^{-|\uz|^2} d\ux\,d\uy\\ &=\frac{1}{\pi^m}\int\limits_{\CC^m}\left[\left(\sum\limits_{|\alpha|=s}a_{\alpha}^{\dagger}(\overline{\uz}^{\CC})^\alpha\right)\left(\sum\limits_{|\beta|=s}a_\beta\uz^\beta\right)\right]_0e^{-|\uz|^2} d\ux\,d\uy\\
&=\frac{1}{\pi^m}\sum\limits_{|\alpha|=s}\sum\limits_{|\beta|=s}\left[a_{\alpha}^{\dagger}a_{\beta}\right]_0\int\limits_{\CC^m}(\overline{\uz}^{\CC})^\alpha\uz^\beta e^{-|\uz|^2} d\ux\,d\uy.
\end{align*}

We solve the integral by transforming the complex coordinates to polar coordinates, i.e. $z_j = r_je^{i\varphi_j}$, $j = 1,\dots,m$. Then,

\begin{align*}
\left\|P_s(\uz)\right\|_{\mathcal{F}^2}^2 =\frac{1}{\pi^m} & \sum\limits_{|\alpha|=s}\sum\limits_{|\beta|=s}\left[a_{\alpha}^{\dagger}a_{\beta}\right]_0\int\limits_{[0,\infty)^m}\int\limits_{[0,2\pi]^m}r_1^{\alpha_1+\beta_1}\dots r_m^{\alpha_m+\beta_m}\\
&\cdot e^{i(\beta_1-\alpha_1)}\dots e^{i(\beta_m-\alpha_m)}e^{-r_1^2-\dots-r_m^2}r_1\dots r_m d\underline{\varphi} d\underline{r}
\end{align*}

The integral $\int\int\dots d\underline{\vphi}d\underline{r}$ is 0 if $\alpha_j\neq\beta_j$ for any $j=1,\dots,m$. So, we get with $\int_0^\infty r^{2n+1}e^{-r^2} dr=\frac{n}{2}$,

\begin{align*}
\left\|P_s(\uz)\right\|_{\mathcal{F}^2}^2 &=\frac{1}{\pi^m}\sum\limits_{|\alpha|=s}\left[a_{\alpha}^{\dagger}a_{\alpha}\right]_0(2\pi)^m\int\limits_{[0,\infty)^m}r_1^{2\alpha_1+1}\dots r_m^{2\alpha_m+1}e^{-r_1^2-\dots-r_m^2} d\underline{r}\\
&=2^m\sum\limits_{|\alpha|=s}|a_{\alpha}|_0^2\prod\limits_{j=1}^m\frac{\alpha_j!}{2}\\
&=\sum\limits_{|\alpha|=s}|a_{\alpha}|_0^2\alpha!
\end{align*}

\item We use the generalization of the Binomial theorem, 

\begin{equation}\label{eq:BinomialTheorem}
|\uz|_0^{2s}=\left(|z_1|^2+\dots+|z_m|^2\right)^s=\sum\limits_{|\alpha|=s}\frac{s!}{\alpha!}|\uz|^{2\alpha},
\end{equation}

and Cauchy-Schwartz (CS) to get

\begin{align*}
\left|P_s(\uz)\right|_0^2 &=\left|\sum\limits_{|\alpha|=s}a_{\alpha}\uz^{\alpha}\right|_0^2\\
&\leq\left(\sum\limits_{|\alpha|=s}\left|a_{\alpha}\uz^{\alpha}\right|_0\right)^2 =\left(\sum\limits_{|\alpha|=s}\sqrt{\frac{\alpha!}{s!}}\left|a_{\alpha}\right|_0\sqrt{\frac{s!}{\alpha!}}\left|\uz^\alpha\right|\right)^2\\
&\stackrel{CS}{\leq}\left(\frac{1}{s!}\sum\limits_{|\alpha|=s}\alpha!\left|a_{\alpha}\right|_0^2\right)\left(\sum\limits_{|\alpha|=s}\frac{s!}{\alpha!}\left|\uz^\alpha\right|^2\right)\\
&\stackrel{(i),(\ref{eq:BinomialTheorem})}{=}\frac{1}{s!}\left\|P_s(\uz)\right\|^2_{\mathcal{F}^2}\left|\uz\right|_0^{2s}. \quad \square
\end{align*}

\end{enumerate}

\begin{lem}\label{lem:sumPsi-infty}

Let $\Psi_{l,k,j}$ be defined as in (\ref{eq:Psi}). Then, 

\[\sup\limits_{\uz\in K}\sum\limits_{l=0}^\infty\sum\limits_{k=0}^\infty\sum\limits_{j=1}^{\mathrm{dim}(M_l^+(k))}|\Psi_{l,k,j}(\uz)|_0^2<\infty\]

for any compact set $K\subset\CC^m$.

\end{lem}

Proof:

Let $\mathrm{SUP}=\sup\limits_{\uz\in K}\sum\limits_{l=0}^\infty\sum\limits_{k=0}^\infty\sum\limits_{j=1}^{\mathrm{dim}(M_l^+(k))}|\Psi_{l,k,j}(\uz)|_0^2$. We first note that each $\Psi_{l,k,j}(\uz)=\frac{1}{\sqrt{\gamma_{l,k}}}\uz^lP_k^{(j)}(\uz)$ is a homogeneous $\C_m^{\CC}$-polyomial of degree $l+k$. Hence, with Lemma \ref{lem:Pk-z} (ii), we get

\[\mathrm{SUP}\leq\sup\limits_{\uz\in K}\sum\limits_{l=0}^\infty\sum\limits_{k=0}^\infty\sum\limits_{j=1}^{\mathrm{dim}(M_l^+(k))}\frac{1}{(l+k)!}\left\|\Psi_{l,k,j}(\uz)\right\|^2_{\mathcal{F}^2}\left|\uz\right|_0^{2l+2k}.\]

We know that $\Psi_{l,k,j}(\uz)=\big(\mathcal{B}\phi_{l,k,j}\big)(\uz)$ (cf. Theorem \ref{th:zlPk} and the proof of Corollary \ref{cor:F2basis}) and that

\[\left\|\mathcal{B}f\right\|_{\mathcal{F}^2}^2=\frac{1}{\sqrt{2\pi}^m}\left\|f\right\|_{L^2}^2\]

for all $f\in L^2$ (cf. Corollary \ref{cor:SegalBargmannIsometry}). Hence,

\[\left\|\Psi_{l,k,j}\right\|^2_{\mathcal{F}^2}=\frac{1}{\sqrt{2\pi}^m}\left\|\phi_{l,k,j}\right\|^2_{L^2}=\frac{1}{\sqrt{2\pi}^m}.\]

We also know that $\mathrm{dim}(M_l^+(k))=\binom{m+k-2}{k}$, cf. section \ref{sub:Pk}. Together, we get

\begin{align*}
\mathrm{SUP}&\leq\frac{1}{\sqrt{2\pi}^m}\sup\limits_{\uz\in K}\sum\limits_{l=0}^\infty\sum\limits_{k=0}^\infty\binom{m+k-2}{k}\frac{1}{(l+k)!}\left|\uz\right|_0^{2l+2k}\\
&\leq\frac{1}{\sqrt{2\pi}^m}\sup\limits_{\uz\in K}\left(\sum\limits_{l=0}^\infty\frac{1}{l!}\left|\uz\right|_0^{2l}\right)\left(\sum\limits_{k=0}^\infty\binom{m+k-2}{k}\frac{1}{k!}\left|\uz\right|_0^{2k}\right)\\
&=\frac{1}{\sqrt{2\pi}^m}\sup\limits_{\uz\in K}\left(\sum\limits_{l=0}^\infty\frac{1}{l!}\left|\uz\right|_0^{2l}\right)\left(\sum\limits_{k=0}^\infty\binom{m+k-2}{k}\frac{1}{2^{mk}k!}\left(2^m\left|\uz\right|_0^2\right)^k\right).\\
\end{align*}

It can be shown via induction that $\binom{m+k-2}{k}\leq 2^{mk}$ for all $k\in\mathbb{N}_0$. Hence,

\begin{align*}
\mathrm{SUP}&\leq\frac{1}{\sqrt{2\pi}^m}\sup\limits_{\uz\in K}\left(\sum\limits_{l=0}^\infty\frac{1}{l!}\left|\uz\right|_0^{2l}\right)\left(\sum\limits_{k=0}^\infty\frac{1}{k!}\left(2^m\left|\uz\right|_0^2\right)^k\right)\\
&=\frac{1}{\sqrt{2\pi}^m}\sup\limits_{\uz\in K}e^{|\uz|_0^2}\cdot e^{2^m\left|\uz\right|_0^2}<\infty. \quad \square
\end{align*}

We are now fully equipped to show convergence of the series expansion (\ref{eq:Psi}).

\begin{prop}\end{prop}\label{prop:convergence}

Let $\phi_{l,k,j}$ be defined as in (\ref{eq:phi2}) and let $\Psi_{l,k,j}$ be defined as in (\ref{eq:Psi}). Then, for each compact set $K\subset\CC^m$,

\[\sup\limits_{z\in K}\left|\sum\limits_{l=0}^\infty\sum\limits_{k=0}^\infty\sum\limits_{j=1}^{\mathrm{dim}(M_l^+(k))}\Psi_{l,k,j}(\uz)\langle\phi_{l,k,j},f\rangle_{L^2}\right|_0<\infty.\]

\end{prop}

Proof:

Let $\mathrm{SUM}=\sum\limits_{l=0}^\infty\sum\limits_{k=0}^\infty\sum\limits_{j=1}^{\mathrm{dim}(M_l^+(k))}\Psi_{l,k,j}(\uz)\langle\phi_{l,k,j},f\rangle_{L^2}$. Since $|\cdot|_0$ is submultiplicative, we have

\begin{align*}
\left|\mathrm{SUM}\right|_0 &\leq\sum\limits_{l=0}^\infty\sum\limits_{k=0}^\infty\sum\limits_{j=1}^{\mathrm{dim}(M_l^+(k))}\left|\Psi_{l,k,j}(\uz)\langle\phi_{l,k,j},f\rangle_{L^2}\right|_0\\
&\leq\sum\limits_{l=0}^\infty\sum\limits_{k=0}^\infty\sum\limits_{j=1}^{\mathrm{dim}(M_l^+(k))}\left|\Psi_{l,k,j}(\uz)\right|_0\cdot\left|\langle\phi_{l,k,j},f\rangle_{L^2}\right|_0
\end{align*}

We now use Proposition \ref{prop:fg0} and $\left\|\phi_{l,k,j}\right\|_{L^2}=1$, so

\begin{align*}
\left|\mathrm{SUM}\right|_0 &\leq\sum\limits_{l=0}^\infty\sum\limits_{k=0}^\infty\sum\limits_{j=1}^{\mathrm{dim}(M_l^+(k))}\left|\Psi_{l,k,j}(\uz)\right|_02^m\left\|\phi_{l,k,j}\right\|_{L^2}\left\|f\right\|_{L^2}\\
&=2^m\left\|f\right\|_{L^2}\sum\limits_{l=0}^\infty\sum\limits_{k=0}^\infty\sum\limits_{j=1}^{\mathrm{dim}(M_l^+(k))}\left|\Psi_{l,k,j}(\uz)\right|_0\\
&\leq 2^m\left\|f\right\|_{L^2}\sqrt{\sum\limits_{l=0}^\infty\sum\limits_{k=0}^\infty\sum\limits_{j=1}^{\mathrm{dim}(M_l^+(k))}\left|\Psi_{l,k,j}(\uz)\right|_0^2}.
\end{align*}

Together with Lemma \ref{lem:sumPsi-infty} the proof is complete. $\square$

Proposition \ref{prop:convergence} shows that $\sum\limits_{l=0}^\infty\sum\limits_{k=0}^\infty\sum\limits_{j=1}^{\mathrm{dim}(M_l^+(k))}\Psi_{l,k,j}(\uz)\langle\phi_{l,k,j},f\rangle_{L^2}$ is absolutely convergent locally uniformly in $\uz\in\CC^m$. Since $\mathcal{B}f$ is the uniform limit of the triple sum on every compact subset of $\CC^m$, it is well-defined and $\mathcal{B}$ can be represented as

\begin{align*}
\big(\mathcal{B}f\big)(\uz) &=\sum\limits_{l=0}^\infty\sum\limits_{k=0}^\infty\sum\limits_{j=1}^{\mathrm{dim}(M_l^+(k))}\Psi_{l,k,j}(\uz)\langle\phi_{l,k,j},f\rangle_{L^2}\\
&=\sum\limits_{l=0}^\infty\sum\limits_{k=0}^\infty\sum\limits_{j=1}^{\mathrm{dim}(M_l^+(k))}\Psi_{l,k,j}(\uz)\int\limits_{\RR^m}\overline{\phi_{l,k,j}(\ux)}f(\ux) d\ux\\
&=\int\limits_{\RR^m}\sum\limits_{l=0}^\infty\sum\limits_{k=0}^\infty\sum\limits_{j=1}^{\mathrm{dim}(M_l^+(k))}\Psi_{l,k,j}(\uz)\overline{\phi_{l,k,j}(\ux)}f(\ux) d\ux\\
&=\left<\overline{\sum\limits_{l=0}^\infty\sum\limits_{k=0}^\infty\sum\limits_{j=1}^{\mathrm{dim}(M_l^+(k))}\Psi_{l,k,j}(\uz)\overline{\phi_{l,k,j}}}\;,f\right>_{L^2}\\
&=\left<\sum\limits_{l=0}^\infty\sum\limits_{k=0}^\infty\sum\limits_{j=1}^{\mathrm{dim}(M_l^+(k))}\phi_{l,k,j}\overline{\Psi_{l,k,j}(\uz)}\;,f\right>_{L^2}.
\end{align*}

Thus,

\[T(\ux,\uz)=\sum\limits_{l=0}^\infty\sum\limits_{k=0}^\infty\sum\limits_{j=1}^{\mathrm{dim}(M_l^+(k))}\phi_{l,k,j}(\ux)\overline{\Psi_{l,k,j}}(\uz) =\frac{1}{\sqrt{2\pi}^m}e^{-\frac{\uz\cdot\uz}{2}+\ux\cdot\uz-\frac{\ux\cdot\ux}{4}}\] 

is the kernel of the Segal-Bargmann transform on $L^2(\RR^m,\C_m^{\RR})$.








%

.



\bibliographystyle{spmpsci}      

\bibliography{BernsteinSchufmann-final}   

\end{document}